\documentclass[10pt,letter]{amsart}
\usepackage{amssymb,amsbsy,amsmath,amsfonts,amssymb}
\usepackage{latexsym,euscript,exscale}

\usepackage{times}

\title{Automatic continuity in homeomorphism groups of compact $2$-manifolds}
\author {Christian Rosendal}
\date {November 2005}
\newcommand {\F}{\mathbb F}

\newcommand{\gr}{{\rm Homeo}(\R^n)}

\newcommand{\grM}{{\rm Homeo}(M)}

\newcommand {\N}{\mathbb N}
\newcommand {\Q}{\mathbb Q}
\newcommand {\R}{\mathbb R}

\newcommand{\eps}{\epsilon}

\newcommand{\isom}{\simeq}

\newcommand{\tom} {\emptyset}

\newcommand{\til}{\rightarrow}

\newcommand {\del}{ \; \big| \;}

\newcommand {\ku} {\mathcal}

\newcommand{\ov}{\overline}
\newcommand{\inv}{^{-1}}

\newtheorem{thm}{Theorem}[section]

\newtheorem{lemme}[thm]{Lemma}

\newtheorem{fact}[thm]{Fact}

\newtheorem{sublemma}[thm]{Sublemma}

\begin{document}
\maketitle

\begin{abstract}
We show that any homomorphism from the homeomorphism group of a compact $2$-manifold, with the compact-open topology, or equivalently, with the topology of uniform convergence, into a separable topological group is automatically continuous.
\end{abstract}

\section{Introduction}

It is well-known and easy to see that for any compact metric space $(X,d)$, its group of homeomorphisms is a separable complete metric group when equipped with the topology of uniform convergence or equivalently with the compact open topology. In fact, a compatible right-invariant metric on ${\rm Homeo}(X,d)$ is given by $d_\infty(g,f)=\sup_{x\in X}d(g(x),f(x))$, and a complete metric by $d_\infty'(g,f)=d_\infty(g,f)+d_\infty(g\inv, f\inv)$. We denote by $B(x,\eps)$ the open ball of radius $\eps$ around $x$ and by $\ov B(x,\eps)$ the corresponding closed ball.

If $g\in {\rm Homeo}(X,d)$, we denote by ${\rm supp}^\circ(g)$ the open set $\{x\in X\del g(x)\neq x\}$ and by ${\rm supp}(g)$ its closure, which we call the {\em support} of $g$. 

We intend to show here that in the case of compact $2$-manifolds, this group topology is intrinsically given by the underlying discrete or abstract group, in the sense that any homomorphism $\pi$ from this group into a separable group is continuous.

\begin{thm}
Let $M$ be a compact $2$-manifold and $\pi:{\rm Homeo}(M)\til H$ a homomorphism into a separable group. Then $\pi$ is automatically continuous when ${\rm Homeo}(M)$ is equipped with the compact-open topology.
\end{thm}

Let us first note the following simple fact, which helps to clear up the situation.

\begin{fact}Suppose $G$ is a topological group. Then the following conditions are equivalent.
\begin{enumerate}
\item Any homomorphism $\pi:G\til {\rm Homeo}([0,1]^\N)$ is continuous,
\item any homomorphism $\pi:G\til H$ into a separable group is continuous.
\end{enumerate}
\end{fact}

\begin{proof}As $[0,1]^\N$ is a compact metric space, its homeomorphism group is a (completely metrisable) separable group in the compact-open topology, so (1) is a special case of (2).

For the other implication, suppose that (1) holds and let $H$ be separable. Let $N$ be the closed normal subgroup of $H$ consisting of all elements that cannot be separated from the identity by an open set and let $H/N$ be the quotient topological group, which is Hausdorff and separable, and, in particular, any non-empty open set covers the group by countably many translates. However, it is an old result (see I.I. Guran \cite{guran}) that for Hausdorff groups this condition is equivalent to being topologically isomorphic to a subgroup of a direct product of separable metric groups, or equivalently, second countable Hausdorff groups (by the Birkhoff-Kakutani metrisation Theorem). Also, a result of Uspenski\u\i  \; \cite{uspenskij} states that any separable metric group is topologically isomorphic to a subgroup of ${\rm Homeo}([0,1]^\N)$, and we can therefore, see $H/N$ as a subgroup of some power of ${\rm Homeo}([0,1]^\N)$. Thus, as a mapping into the Tikhonov product is continuous if and only if the composition with each coordinate projection is continuous, $\pi$ composed with the quotient mapping is continuous, and hence by the choice of $N$, also $\pi$ is continuous.\end{proof}

However, we shall not use this result in any way, but instead simplify matters by not be working with arbitrary homomorphisms, but rather with arbitrary subsets of the group satisfying a certain algebraic largeness condition. Let $G$ be a group and $W\subseteq G$ be a symmetric set. We say that $W$ is {\em countably syndetic} if there are countably many left-translates of $W$ whose union cover $G$. Moreover if $G$ is a topological group, we say that $G$ is {\em Steinhaus} if for some $k\geq 1$ and all symmetric, countably syndetic $W\subseteq G$, $\rm Int(W^k)\neq \tom$. It is not hard to prove (see, e.g., \cite{rossol}) that Steinhaus groups satisfy the equivalent conditions of the above fact, and this is the condition that we will verify. Note however the order of quantification; the $k$ is universal for all symmetric, countably syndetic $W$. Indeed, the group ${\rm Homeo}_+(S^1)$ equipped with the trivial topology $\tau=\{\tom,{\rm Homeo}_+(S^1)\}$ satisfies the condition when we have inversed the quantifiers, but the identity homomorphism into itself equipped with the compact-open topology is obviously discontinuous.

It is instructive also to consider from which groups one can construct discontinuous homeomorphisms. Of course the first case that comes to mind is $(\R,+)$, on which one can with the help of a Hamel basis, i.e., a basis for $\R$ as a $\Q$-vector space, construct discontinuous automorphisms, and, in fact, construct group isomorphisms between $\R$ and $\R^2$.

\section{The proof}

\subsection{Commutators}
We shall first prove a general lemma about homeomorphisms of $\R^n$.

\begin{lemme}\label{commutator}Suppose that  $g\in {\rm Homeo}(\R^n)$ has compact support. Then there are $f,h\in \gr$ with compact support such that $g=[f,h]=fhf\inv h\inv$.
\end{lemme}

\begin{proof} Fix some open ball $U_0\subseteq \R^n$ containing the support of $g$ and let $(U_m)$ be a sequence of disjoint open balls such that for some distinct $x_0$ and $x_1$ in $\R^n$, the sequences $(\ov U_m)_{m\geq 0}$ and $(\ov U_{-m})_{m\geq 0}$ converge in the Vietoris topology to $x_0$ and $x_1$ respectively. We can now find a shift $h\in \gr$ with compact support, i.e., such that $h[U_m]=U_{m+1}$ and define our $f$ by letting $f|U_m=h^mgh^{-m}|U_m$ for $m\geq 0$ and setting $f={\rm id}$ everywhere else. We now see that for $m>0$,
\begin{displaymath}
hf\inv h\inv |U_m=h(h^{m-1}g\inv h^{-m+1})h\inv|U_m= h^mg\inv h^{-m}|U_m,
\end{displaymath}
and for $m\leq 0$,
\begin{displaymath}
hf\inv h\inv|U_m =h\,{\rm id}\,h\inv|U_m={\rm id}|U_m,
\end{displaymath}
while $hf\inv h\inv={\rm id}$ everywhere else. Therefore, $f\cdot hf\inv h\inv |U_m={\rm id}|U_m$ for $m>0$, $f\cdot hf\inv h\inv |U_0=f|U_0=g|U_0$, $f\cdot hf\inv h\inv |U_m={\rm id}|U_m$ for $m<0$, and $fhf\inv h\inv ={\rm id}$ everyhere else. This shows that $g=[f,h]=fhf\inv h\inv$. 
\end{proof}

We notice that in the proof above we used $f$ and $h$ with slightly bigger support than $g$. I believe it is an open problem whether this can be avoided and indeed it seems to be a much harder problem. We can restate the problem as follows. Can every homeomorphism of $[0,1]^n$ that fixes the boundary pointwise be written as a commutator of $f$ and $h$ that also fixes the boundary pointwise? What happes if we replace pointwise by setwise? Let us mention that the first question has a positive answer in dimension $1$ as, for example, the group of orientation preserving homeomorphisms of $[0,1]$ has a comeagre conjugacy class \cite{kecros}. 
The above result slightly strengthens a result of Mather \cite{mather} saying that the homology groups of the group of homeomorphisms $\R^n$ with compact support vanish. One can of course also extend the lemma to $[0,\infty[\times \R^{n-1}$ and thus also improve the result of Rybicki \cite{rybicki}.

\subsection{Countably syndetic sets}
We will now prove some properties of coutably syndetic sets in the homeomorphism groups of arbitrary manifolds. These results will allow us to completely solve our problem for compact two-dimensional manifolds and provide techniques for higher dimensions. So let $M$ be a manifold of dimension $n$ and fix a compatible complete metric $d$ on $M$.

In the following we fix a countably syndetic symmetric subset $W\subseteq \grM$ and a sequence $k_m\in \grM$ such that $\bigcup_mk_mW=\grM$.

\begin{lemme}\label{patched}
For all distinct $y_1,\ldots,y_p\in M$ and $\eps>0$, there are $\eps>\delta>0$ and $z_i\in B(y_i,\eps)$ such that if $g\in \grM$ has support contained in $D=\bigcup_{i=1}^p \ov B(z_i,\delta)$, then $g\in W^{16}$.
\end{lemme}

\begin{proof}We notice that it is enough to find $z_i\in B(y_i,\eps)$ and open neighbourhoods $U_i$ of $z_i$ such that if $g\in \grM$ has support contained in $\bigcup_iU_i$, then $g\in W^{16}$.  We choose some open neighbourhood of $y_i$, $E_i\subseteq B(y_i,\eps)$, that is homeomorphic to $]0,\eps[^2$. We also suppose that the sets $E_i$ are $4\eps$-separated.
We will also temporarily transport the standard euclidian metric from $]0,\eps[^2$ to each of the sets $E_i$. As we will be working separately on each of $E_i$, this will not cause a problem. Thus in the following, the notation $B(x,\beta)$ will refer to the balls in the transported euclidian metric, which we denote by $d$.

\begin{sublemma}\label{semi-patched}
For all $u_i\in E_i$ and $\gamma>0$ such that $d(u_i,\partial E_i)>2\gamma$, there are $\gamma>\alpha>0$ and $x_i\in \partial B(u_i,\gamma)$ such that if $g\in \grM$ has support contained in $A=\bigcup_{i=1}^p \ov B(x_i,\alpha)\cap\ov B(u_i,\gamma)$, then there is an $h\in W^2$ with support contained in $\bigcup_{i=1}^p\ov B(u_i,\gamma)$ such that $g|A=h|A$.
\end{sublemma}

\begin{proof}Let $u_1,\ldots,u_p$ be given. We fix for each $i\leq p$ a sequence of distinct points $x^i_m\in \partial B(u_i,\gamma)$ converging to some point $x^i_\infty\in\partial B(u_i,\gamma)$ and choose a sequence $\frac{\gamma}{2}>\alpha_m>0$ such that $B(x^i_m,\alpha_m)\cap B(x^i_l,\alpha_l)=\tom$ for any $m\neq l$. Thus, as $\alpha_m\til 0$, we have that if $g_m\in \grM$ has support only in 
$$
A_m=\big(\ov B(x^1_m,\alpha_m)\cap \ov B(u_1,\gamma)\big)\cup\ldots \cup \big(\ov B(x^p_m,\alpha_m)\cap \ov B(u_p,\gamma)\big)
$$ 
for each $m\geq 0$, then there is a homeomorphism $g\in \grM$, whose support is contained in $C=\ov B(u_1,\gamma)\cup\ldots\cup \ov B(u_p,\gamma)$, such that $g| A_m=g_m|A_m$. We claim that for some $m_0\geq 0$, if $g\in \grM$ has support contained in $A_{m_0}$, then there is an element $h\in k_{m_0}W$, with support contained in $C$, such that $g|A_{m_0}=h|A_{m_0}$. Assume toward a contradiction that this is not the case. Then for every $m$ we can find some $g_m\in \grM$ with support contained in  $A_m$ such that for all $h\in k_mW$, if ${\rm supp}(h)\subseteq C$, then $g_m|A_m\neq h|A_m$. But then letting $g\in \grM$ have support in $C$ and agree with each $g_m$ on $A_m$ for each $m$, we see that if $h\in k_mW$ has support in $C$, then $g$ disagrees with $h$ on  $A_m$. Therefore, $g$ cannot belong to any $k_mW$, contradicting that these cover $\grM$. Suppose that $m_0$ has been chosen as above and denote $x^i_{m_0}$ by $x_i$, $A_{m_0}$ by $A$, and $\alpha_{m_0}$ by $\alpha$.

Then for any $g\in \grM$ with support contained in $A$, there is an element $h\in W^2$ with support contained in $C$ such that $g|A=h|A$ for all $i\leq p$. To see this, it is enough to notice that we can find $h_0,h_1\in k_{m_0}W$, with ${\rm supp}(h_0),{\rm supp}(h_1)\subseteq C$, such that  $g|A=h_1|A$ and ${\rm id}|A=h_0|A$. But then $h_0\inv h_1\in (k_{m_0}W)\inv k_{m_0}W=W\inv W=W^2$ and $g|A={\rm id}\,g|A=h_0\inv h_1|A$.
\end{proof}

We will first apply Sublemma \ref{semi-patched} to the situation where $u_i=y_i$ and $\gamma>0$ is sufficiently small. We thus obtain $\gamma>\alpha>0$ and $x_i\in \partial B(y_i,\gamma)$ such that if $g\in \grM$ has support contained in $A=\bigcup_{i=1}^p \ov B(x_i,\alpha)\cap\ov B(y_i,\gamma)$, then there is an $h\in W^2$ with support contained in $\bigcup_{i=1}^p\ov B(y_i,\gamma)$ such that $g|A=h|A$. 

Now pick $y'_i\in B(x_i,\alpha)\cap B(y_i,\gamma)$ and $\gamma'>0$ such that $B(y'_i,2\gamma')\subseteq B(x_i,\alpha)\cap B(y_i,\gamma)$. We now apply Lemma \ref{semi-patched} once again to this new situation, in order to obtain $\gamma'>\alpha'>0$ and $x'_i\in \partial B(y'_i,\gamma')$ such that if $g\in \grM$ has support contained in $A'=\bigcup_{i=1}^p \ov B(x'_i,\alpha')\cap\ov B(y'_i,\gamma')$, then there is an $h\in W^2$ with support contained in $\bigcup_{i=1}^p\ov B(y'_i,\gamma')$ such that $g|A'=h|A'$. 

Now clearly there is a homeomorphism $a\in \grM$ whose support is contained in $A=\bigcup_{i=1}^p \ov B(x_i,\alpha)\cap\ov B(y_i,\gamma)$ such that $a[B(y'_i,\gamma')]=B(x'_i,\alpha')$ and 
$$
a[\ov B(y'_i,\gamma')\cap \ov B(x'_i,\alpha')]=\ov B(y'_i,\gamma')\cap \ov B(x'_i,\alpha'),
$$ 
and hence we can also find such an $a$ in $W^2$, except that its support may now be all of $\bigcup_{i=1}^p\ov B(y_i,\gamma)$.

We therefore have that if $g\in \grM$ has support contained in $A'$, then $a\inv ga$ also has support contained in $A'$, and so there is an $h\in W^2$ with support contained in $\bigcup_{i=1}^p\ov B(y'_i,\gamma')$ such that $a\inv ga|A'=h|A'$. But then $g|A'=aha\inv |A'$, while 
$$
{\rm supp}(aha\inv)=a[{\rm supp}(h)]\subseteq a[\bigcup_{i=1}^p\ov B(y'_i,\gamma')]=\bigcup_{i=1}^p\ov B(x'_i,\alpha').
$$ 

We now notice that $aha\inv \in W^6$, and thus that if $g\in \grM$ has support contained in $A'=\bigcup_{i=1}^p \ov B(x'_i,\alpha')\cap\ov B(y'_i,\gamma')$, then there is some $f\in W^6$ with support contained in $\bigcup_{i=1}^p\ov B(x'_i,\alpha')$ such that $g|A'=f|A'$.

Now suppose finally that $g\in \grM$ is any homeomorphism having support contained in $\bigcup_{i=1}^p B(x'_i,\alpha')\cap B(y'_i,\gamma')$. Since the sets $B(x'_i,\alpha')\cap B(y'_i,\gamma')$ are homeomorphic to $\R^n$, working separately on each of these sets and noticing that $g$ has compact support, we can invoke Lemma \ref{commutator} to write $g$ as a commutator $[b,c]$ for some $b,c\in \grM$ whose supports are contained in $\bigcup_{i=1}^p B(x'_i,\alpha')\cap B(y'_i,\gamma')\subseteq A'$. Find now $h\in W^2$ agreeing with $b$ on $A'$ and with support contained in $\bigcup_{i=1}^p\ov B(y'_i,\gamma')$, and, similarly, find $f\in W^6$ agreeing with $c$ on $A'$ and with support contained in $\bigcup_{i=1}^p\ov B(x'_i,\alpha')$. Then the set of common support of $h$ and $f$ is included in $A'$ on which they agree with $b$ and $c$ respectively, and we have therefore that $[h,f]=hfh\inv f\inv=bcb\inv c\inv =g$. In other words, $g\in W^{16}$. We can therefore finish the proof by choosing some $z_i\in B(x'_i,\alpha')\cap B(y'_i,\gamma')$ and letting $U_i=B(x'_i,\alpha')\cap B(y'_i,\gamma')$.
\end{proof}

\subsection{Circular orders}
In order to simplify notation, we will consider {\em circular orders} on finite sets. Since we are really just interested in simplifying notation, let me just say what a circular order is in terms of an example, namely, $S^1$. For $x,y,z$ distinct points on $S^1$,  $y$ is said to be between $x$ and $z$, in symbols $B(x,y,z)$, if going counterclockwise around $S^1$ from $x$ to $y$, one does not pass through $z$. Thus a circular order is just a circular betweeness relation. When $B$ is a circular order on a finite set $\F$, we denote for each $x\in \F$ its immediate successor and immediate predecessor, i.e., the first elements encountered by going respectively counterclockwise and clockwise around $\F$, by $x^+$ and $x^-$. So, e.g., $(x^+)^-=x$.

\subsection{A quantitative annulus theorem}
Fix three points $v_0,v_1,v_2\in \R^2$ such that for $i\neq j$, $d(v_i,v_j)=1$,  and denote by $\triangle$ the $2$-cell consisting of the points lying within the triangle $\triangle v_0v_1v_2$. Suppose also that the barycenter of $\triangle$ lies at the origin, so that for all $\lambda>0$, $\lambda \triangle$ and $\triangle$ are concentric triangles, the former with sidelengths $\lambda$.

\begin{lemme}
Let $\phi:(1-2\eta)\triangle\til \triangle$ be a homeomorphic embedding satisfying 
$$
\sup_{x\in(1-2\eta)\triangle}d(x,\phi(x))<\frac\eta{100},
$$
where $\eta<\frac1{1000}$. Then there is a homeomorphism $\psi:\triangle\til \triangle$ that is the identity outside of $(1-\eta)\triangle$, with $\sup_{x\in\triangle}d(x,\psi(x))<100\eta$, and such that $\psi\circ\phi|_{(1-2\eta)\triangle}={\rm id}$.
\end{lemme}

\begin{proof}Let $\partial(1-\eta)\triangle$ be the boundary of $(1-\eta)\triangle$ and pick a finite set of points $\F$ containing $(1-\eta)v_0,(1-\eta)v_1,(1-\eta)v_2$ and lying in $\partial(1-\eta)\triangle$, such that when $\F$ is equipped with the circular order obtained from going counterclockwise around $\partial(1-\eta)\triangle$, we have $d(x,x^+)\in ]20\eta,21\eta[$ for all $x\in \F$. As $\triangle$ is equilateral, $d(x,y)>20\eta$ for all $x\neq y$ in $\F$. 

Let now $C=\phi[\partial(1-2\eta)\triangle]$ be the image of the boundary of $(1-2\eta)\triangle$, so $C$ is a simple closed curve. Choose also for each $x\in \F$ a point $\hat x\in C$ such that the distance $d(x,\hat x)$ is minimal. Since $\sup_{x\in (1-2\eta)\triangle} d(x,\phi(x))<\frac\eta{100}$ and 
$$
\frac \eta3<d(x,\partial(1-2\eta)\triangle)<\frac{2\eta}3
$$
for all $x\in \partial(1-\eta)\triangle$, also $d(x,\hat x)<\eta$ and $d(C,\partial(1-\eta)\triangle)>\frac\eta4$.

For all $x\in \F$, denote by $\alpha_x$ the straight (oriented) line segment from $x$ to $\hat x$ and by $\beta_x$ the straight line segment from $x$ to $x^+$. We also let $\gamma'_x$ be the shortest path in $\partial(1-2\eta)\triangle$ from $\phi\inv (\hat x)$ to $\phi\inv (\widehat{x^+})$ and put $\gamma_x=\phi[\gamma'_x]$. 

By definition of $\hat x$, $\alpha_x$ intersects $C$ exactly in $\hat x$, intersects $\partial(1-\eta)\triangle$ in exactly $x$, and therefore $\alpha_x$ and $\gamma_y$ intersect only if $y=x^-$ or $y=x$. Similarly, none of the paths $\beta_x$ and $\gamma_y$ intersect as they lie in $\partial(1-\eta)\triangle$ and $C$ respectively. Therefore, for any $x\in \F$, $\ku C_x=\alpha_x\centerdot\gamma_x\centerdot\bar\alpha_{x^+}\centerdot\bar\beta_x$ is a simple closed curve beginning and ending at $x$. Here $\bar\alpha$ denotes the reverse path of $\alpha$ and $\centerdot$ the concatenation of paths. By the Sch\"onflies Theorem, $\R^2\setminus \ku C_x$ has exactly two components, one unbounded and the other $U_x$ bounded, homeomorphic with $\R^2$ and with boundary $\ku C_x$. Moreover, as the diameter of $\ku C_x$ is bounded by $30\eta$, $\ku C_x$ intersects $\partial(1-\eta)\triangle$ in exactly $\beta_x$, and the diameter of $\partial(1-\eta)\triangle\setminus\beta_x$ is $1-\eta>30\eta$, this means that $\partial(1-\eta)\triangle\setminus\beta_x$ lies in the unbounded component. Therefore, if $R_x=\overline U_x=U_x\cup \ku C_x$, we have for $x\neq y$
\begin{displaymath}
\begin{split}
R_x\cap R_y =\left \{
\begin{array}{ll}
\tom \;\;\;\;\;\;\;\;\;\;\;\;&\textrm{ if } y\neq x^+ \textrm { and } y\neq x^-\\
\alpha_y\;\;\;\;\;\;\;\;\;\;\;\;&\textrm{ if } y= x^+\\
\alpha_x\;\;\;\;\;\;\;\;\;\;\;\;&\textrm{ if } y= x^-
\\\end{array}\right. 
\end{split}
\end{displaymath}
We can now define $\psi:\triangle\til \triangle$ by letting $\psi=\phi\inv$ on $\phi[(1-2\eta)\triangle]$, $\psi={\rm id}$ on $\triangle\setminus (1-\eta)\triangle$, and, moreover, along the boundaries of $R_x$ construct $\psi$ as follows: $\psi[\alpha_x]$ is the straight line segment from $x$ to $\phi\inv (\hat x)$, $\psi[\gamma_x]=\gamma_x'$, and $\psi[\beta_x]=\beta_x$. Then 
$$
\psi[\ku C_x]=\psi[\alpha_x\centerdot\gamma_x\centerdot\bar\alpha_{x^+}\centerdot\bar\beta_x]=
\psi[\alpha_x]\centerdot\psi[\gamma_x]\centerdot\overline{\psi[\alpha_{x^+}]}\centerdot\overline{\psi[\beta_x]}=
\psi[\alpha_x]\centerdot\gamma'_x\centerdot\overline{\psi[\alpha_{x^+}]}\centerdot\overline\beta_x
$$
is the boundary of a region $K_x$ homeomorphic to the unit disk $D^2$ and hence, by Alexander's Lemma, the homeomorphism $\psi$ from $\ku C_x=\alpha_x\centerdot\gamma_x\centerdot\bar\alpha_{\hat x}\centerdot\bar\beta_x$ to $\psi[\alpha_x]\centerdot\gamma'_x\centerdot\overline{\psi[\alpha_{x^+}]}\centerdot\overline\beta_x$ extends to the regions that they bound, i.e., to a homeomorphism of $R_x$ with $K_x$. This finishes the description of $\psi$ and it therefore only remains to see that $\sup_{x\in\triangle}d(x,\psi(x))<100\eta$. Since $\psi=\phi\inv$ on $\phi[(1-2\eta)\triangle]$ and $\psi={\rm id}$ on $\triangle\setminus(1-\eta)\triangle$ it is enough to consider what $\psi$ does to $x\in (1-\eta)\triangle\setminus\phi[(1-2\eta)\triangle]\subseteq \bigcup_{x\in \F}R_x$. Now, $\psi[R_x]=K_x$ for all $x\in \F$, and hence it is enough to show that no points in $R_x$ and in $K_x$ are more than $100\eta$ apart. But ${\rm diam}(R_x)<30\eta$ and ${\rm diam}(K_x)<40\eta$, while $R_x\cap K_x\neq \tom$, which gives the desired result. This finishes the proof. 
\end{proof}

\subsection{Patching along a triangulation of a compact $2$-manifold}
As $\grM$ is a separable complete metric group it is not covered by countably many nowhere dense sets (this is the Baire category theorem) and hence $W$ must be dense in some non-empty open set, whereby $W\inv W=W^2$ is dense in some neighbourhood of the identity in $\grM$. So fix some $\eta_1>0$ such that $W^2$ is dense in 
\begin{equation}
V_{\eta_1}=\{g\in \grM \del d_\infty(g,{\rm id})<\eta_1\}.
\end{equation}

It is a well-known fact, first proved rigorously by Tibor Rado, that any compact $2$-manifold can be triangulated. So from now on, we assume that $M$ is a fixed compact $2$-manifold and we pick a triangulation $\{T_1,\ldots, T_m\}$ of $M$ with corresponding homeomorphisms $\chi_i:\triangle\til T_i$. By further triangulating each $T_i$, we can suppose that the diameter of $T_i$ is less than $\frac{\eta_1}{10}$ for all $i$. Moreover, by first modifying the $\chi_i$ along each edge of $\triangle$ and then extending to the interior of $\triangle$ by Alexander's Lemma, we can suppose that the following holds. If $T_i=\chi_i[\triangle]$ and $T_j=\chi_j[\triangle]$ have an edge in common, then $\chi_i$ and $\chi_j$ agree along this edge, i.e., if $\chi_i(v_a)=\chi_j(v_\alpha)$ and $\chi_i(v_b)=\chi_j(v_\beta)$, then for all $t\in [0,1]$, $\chi_i(tv_a+(1-t)v_b)=\chi_j(tv_\alpha+(1-t)v_\beta)$.

\begin{lemme}\label{triangle-support}For all $0<\eta<1$, if $h\in {\rm Homeo}(M)$ has support contained in $$\bigcup_{i=1}^m\chi_i[(1-\eta)\triangle],$$ then $h\in W^{20}$.
\end{lemme}

\begin{proof}
Let $y_i=\chi_i(\vec 0)$ and choose $\eps>0$ such that $\overline B(y_i,\eps)\subseteq \chi_i[(1-\eta)\triangle]$ for all $i\leq m$. By Lemma \ref{patched}, we can find some $0<\delta<\eps$ and $z_i\in B(y_i,\eps)$ such that if $g\in {\rm Homeo}(M)$ has support contained in $\bigcup_{i=1}^m\overline B(z_i,\delta)$ then $g\in W^{16}$.

As $W^2$ is dense in $V_{\eta_1}$, we can find an $f\in W^2$ such that for every $i\leq m$, $f[\chi_i[(1-\eta)\triangle]]\subseteq \overline B(z_i,\delta)$ and thus if $h$ is given as in the statement of the lemma, ${\rm supp}(fhf\inv)=f[{\rm supp}(h)]\subseteq \bigcup_{i=1}^m\overline B(z_i,\eps)$ and thus $g=fhf\inv\in W^{16}$, whence $h\in W^{20}$.
\end{proof}

\begin{lemme}\label{fixing triangles}
Let $\delta,\eta>0$, $\eta<\frac1{1000}$ be such that for $i\leq m$ and $x,y\in \triangle$,
$$d(x,y)<100\eta\til d(\chi_i(x),\chi_i(y))<\delta.$$
Then there is an $\alpha>0$ such that for all $g\in V_\alpha$ there is $\psi\in V_{\delta}\cap W^{20}$ whose support is contained in 
$\bigcup_{i=1}^m\chi_i[ {(1-\eta)\triangle}]$ and such that for all $i\leq m$,
$$\psi\circ g|_{\chi_i[(1-2\eta)\triangle]}={\rm id}.$$
\end{lemme}

\begin{proof}
Fix $\delta$ and $\eta$ as in the Lemma. Then for any  continuous $\phi:\triangle \til \triangle$ such that $\sup_{x\in \triangle} d(x,\phi(x))<100\eta$, we have for every $i\leq m$,
$$\sup_{y\in T_i}d(y,\chi_i\circ\phi\circ\chi_i\inv(y))=\sup_{x\in \triangle}d(\chi_i(x),\chi_i\circ\phi(x))<\delta.$$
Now pick some $\alpha>0$ such that for $g\in V_\alpha$ and $i\leq m$, we have 
$$g\circ\chi_i[(1-2\eta)\triangle]\subseteq \chi_i[\triangle]=T_i,$$
whereby $\chi_i\inv\circ g\circ\chi_i:(1-2\eta)\triangle\til \triangle$, and such that
$$
\sup_{x\in (1-2\eta)\triangle}d(x,\chi_i\inv\circ g\circ\chi_i(x))<\frac\eta{100}.
$$

By the quantitative annulus theorem we can therefore find some homeomorphism $\psi_i:\triangle\til \triangle$ that is the identity outside of $(1-\eta)\triangle$, satisfies $\sup_{x\in \triangle}d(x,\psi_i(x))<100\eta$, and 
$$
\psi_i\circ\chi_i\inv\circ g\circ \chi_i|_{(1-2\eta)\triangle}={\rm id}.
$$
This implies that for each $i\leq m$, $\chi_i\circ \psi_i\circ\chi_i\inv:T_i\til T_i$ is a homeomorphism that is the identity outside of $\chi_i[(1-\eta)\triangle]$, $\sup_{x\in T_i}d(x,\chi_i\circ \psi_i\circ\chi_i\inv(x))<\delta$, and 
$$
\chi_i\circ \psi_i\circ\chi_i\inv\circ g|_{\chi_i[{(1-2\eta)\triangle}]}={\rm id}.
$$
We can therefore define $\psi=\bigcup_{i=1}^m  \chi_i\circ \psi_i\circ\chi_i\inv\in {\rm Homeo}(M)$ and notice that $\psi\in V_\delta$ and $\psi\circ g|_{\chi_i[ {(1-2\eta)\triangle}]}={\rm id}$ for every $i\leq m$. We see that $\psi$ has its support contained within the set $\bigcup_{i=1}^m\chi_i[ {(1-\eta)\triangle}]$ and thus, by Lemma \ref{triangle-support}, $\psi$ belongs to $W^{20}$.
\end{proof}

Fix some $0<\tau<\frac{1}{100}$. We now define the following set of points in $\triangle$: For distinct $i,j=0,1,2$, we put $w_{ij}=(1-10\tau)v_i+10\tau v_j$, $w_{ij}^+=(1-9\tau)v_i+9\tau v_j$, $u_{ij}=(1-\tau)w_{ij}$ and $u_{ij}^+=(1-\tau)w^+_{ij}$. So $w_{ij},w_{ij}^+\in \partial\triangle$, while $u_{ij},u_{ij}^+\in\partial(1-\tau)\triangle$. We also define a number of paths as follows:
\begin{itemize}
\item $\alpha_{ij}$ is the straight line segment from $u_{ij}$ to $w_{ij}$.
\item $\beta_{ij}$ is the straight line segment from $w_{ij}$ to $w^+_{ij}$.
\item $\gamma_{ij}$ is the straight line segment from $u^+_{ij}$ to $w^+_{ij}$.
\item $\zeta_{ij}$ is the straight line segment from $u_{ij}$ to $u^+_{ij}$.
\item $\kappa_{ij}$ is the straight path from $w_{ij}$ to $w_{ji}$.
\item $\omega_{ij}$ is the straight path from $u_{ij}$ to $u_{ji}$.

\item $\xi_{0}$ is the shortest path in $\partial(1-\tau)\triangle$ from $u^+_{02}$ to $u^+_{01}$.
\item $\xi_{1}$ is the shortest path in $\partial(1-\tau)\triangle$ from $u^+_{10}$ to $u^+_{12}$.
\item $\xi_{2}$ is the shortest path in $\partial(1-\tau)\triangle$ from $u^+_{21}$ to $u^+_{20}$.

\item $\theta_{0}$ is the shortest path in $\partial\triangle$ from $w^+_{02}$ to $w^+_{01}$.
\item $\theta_{1}$ is the shortest path in $\partial\triangle$ from $w^+_{10}$ to $w^+_{12}$.
\item $\theta_{2}$ is the shortest path in $\partial\triangle$ from $w^+_{21}$ to $w^+_{20}$.
\end{itemize}
We thus see that 
$$
\ku C_{ij}=\kappa_{ij}\centerdot\overline \alpha_{ji}\centerdot \omega_{ji}\centerdot\alpha_{ij}
$$ 
is a simple closed curve bounding a closed region $R_{ij}=R_{ji}\subseteq \triangle$, 
$$
\ku C^+_{ij}=\overline\beta_{ij}\centerdot\kappa_{ij}\centerdot \beta_{ji}\centerdot\overline\gamma_{ji}\centerdot\overline\zeta_{ji}\centerdot\omega_{ji}\centerdot\zeta_{ij}\centerdot\gamma_{ij}
$$ 
is a simple closed curve bounding a closed region $R^+_{ij}=R^+_{ji}\subseteq \triangle$ that contains $R_{ij}$. 

Notice however that the preceding definitions depend on the choice of $\tau$, which is therefore also the case for the following lemma.

\begin{lemme}\label{band-support}
If $\phi\in {\rm Homeo}(M)$ has support contained in $\bigcup_{l=1}^m\bigcup_{0\leq i<j\leq 2}\chi_l[R_{ij}^+]$, then $\phi\in W^{20}$.
\end{lemme}

\begin{proof}
We notice that for distinct $l,l'$, $\chi_l[R^+_{ab}]\cap \chi_{l'}[R^+_{a'b'}]\neq \tom$ if and only if the triangles $T_l$ and $T_{l'}$ have the edge $\chi_l[\overline{v_av_b}]=\chi_{l'}[\overline{v_{a'}v_{b'}}]$ in common. Moreover, in this case, the set $\chi_l[R^+_{ab}]\cup \chi_{l'}[R^+_{a'b'}]$ is homeomorphic to the unit disk $D^2$ and is contained in an open set homeomorphic to $\R^2$.

So let $A_1,\ldots,A_{\frac{3m}{2}}$ be an enumeration of all the closed sets $\chi_l[R^+_{ab}]\cup \chi_{l'}[R^+_{a'b'}]$ with $\chi_l[R^+_{ab}]$ and $\chi_{l'}[R^+_{a'b'}]$ overlapping and let $U_i\subseteq M$ be an open set containing $A_i$, homeomorphic to $\R^2$.
We can suppose that the $U_i$ are all pairwise disjoint. Moreover, as the diameter of each $T_j$ is at most $\frac{\eta_1}{10}$, the diameter of each $A_i$ is at most $\frac{\eta_1}5$.

The proof is now very much the same as the proof of Lemma \ref{triangle-support}. Let $y_i\in A_i$ and choose $0<\eps<\frac{\eta_1}{5}$ such that $\overline B(y_i,\eps)\subseteq U_i$ for all $i\leq m$. By Lemma \ref{patched}, we can find some $0<\delta<\eps$ and $z_i\in B(y_i,\eps)$ such that if $g\in {\rm Homeo}(M)$ has support contained in $\bigcup_{i=1}^m\overline B(z_i,\delta)$ then $g\in W^{16}$.

As $W^2$ is dense in $V_{\eta_1}$, we can find an $f\in W^2$ such that for every $i\leq \frac{3m}2$, $f[A_i]\subseteq \overline B(z_i,\delta)$ and thus if $\phi$ is given as in the statement of the lemma, 
$$
{\rm supp}(f\phi f\inv)=f[{\rm supp}(\phi)]\subseteq \bigcup_{i=1}^m\overline B(z_i,\eps),
$$ 
and thus $g=f\phi f\inv\in W^{16}$, whence $\phi\in W^{20}$.
\end{proof}

\begin{lemme}\label{vertices left}
There is a $\nu>0$ such that if $g\in V_\nu$ and $g$ is the identity on $\bigcup_{i=1}^m\chi_i[ {(1-\tau)\triangle}]$, then there is a $\phi\in W^{20}$ such that $\phi\circ g$ is the identity on 
$$
\bigcup_{i=1}^m\chi_i[ {(1-\tau)\triangle}]\cup \bigcup_{l=1}^m\;\bigcup_{0\leq i<j\leq 2}\chi_l[R_{ij}].
$$ 
\end{lemme}

\begin{proof}
Consider the closed set $M_0=M\setminus{\rm Int}(\bigcup_{i=1}^m\chi_i[ {(1-\tau)\triangle}])$ and the closed subgroup $H=\{g\in {\rm Homeo}(M)\del g|_{ \bigcup_{i=1}^m\chi_i[ {(1-\tau)\triangle}]}={\rm id}\}$.
Assume that $T_l$ and $T_{l'}$ have an edge in common, i.e., $\chi_l(v_a)=\chi_{l'}(v_{a'})$ and $\chi_l(v_b)=\chi_{l'}(v_{b'})$ for some $a,a',b,b'$. Then $\chi_l[R_{ab}]\cup \chi_{l'}[R_{a'b'}]\subseteq {\rm Int}_{M_0}(\chi_l[R^+_{ab}]\cup \chi_{l'}[R^+_{a'b'}])$. Therefore, we can find some $\nu>0$, not depending on the particular choice of $l,l',a,a',b,b'$, such that for all such choices of $l,l',a,a',b,b'$ and $g\in V_\nu\cap H$ we have
\begin{align}\label{g}
g[\chi_l[R_{ab}]\cup \chi_{l'}[R_{a'b'}]]\subseteq {\rm Int}_{M_0}(\chi_l[R^+_{ab}]\cup \chi_{l'}[R^+_{a'b'}]).
\end{align}
Fix some $g\in V_\nu\cap H$.

Assume now that $\chi_l[\triangle]$ and $\chi_k[\triangle]$ have an edge in common. For concreteness we can suppose that, e.g., $\chi_l(v_0)=\chi_k(v_1)$ and $\chi_l(v_1)=\chi_k(v_2)$. As the covering mappings $\chi_i$ were supposed to agree along their edges, this implies that $\chi_l[\beta_{01}]=\chi_k[\beta_{12}]$, $\chi_l[\kappa_{01}]=\chi_k[\kappa_{12}]$, and $\chi_l[\beta_{10}]=\chi_k[\beta_{21}]$.
Also, as $g\in H$, $g$ is the identity on the paths $\chi_l[\zeta_{01}],\chi_l[\omega_{01}],\chi_l[\zeta_{10}], \chi_k[\zeta_{12}],\chi_k[\omega_{12}]$ and $\chi_k[\zeta_{21}]$.

By consequence, $\chi_l[\zeta_{01}]\centerdot \chi_l[\gamma_{01}]\centerdot\chi_k[\overline\gamma_{12}]\centerdot\chi_k[\overline\zeta_{12}]$ and $\chi_l[\alpha_{01}]\centerdot \chi_k[\overline\alpha_{12}]$ are paths from $\chi_l(u_{01})$ to $\chi_k(u_{12})$ only intersecting in their endpoints. Similarly, $\chi_l[\zeta_{10}]\centerdot \chi_l[\gamma_{10}]\centerdot\chi_k[\overline\gamma_{21}]\centerdot\chi_k[\overline\zeta_{21}]$ and $\chi_l[\alpha_{10}]\centerdot\chi_k[\overline\alpha_{21}]$ are paths from $\chi_l(u_{10})$ to $\chi_k(u_{21})$ only intersecting in their endpoints. This shows that 
$$
\ku K=\chi_l[\zeta_{01}]\centerdot\chi_l[\gamma_{01}]\centerdot\chi_k[\overline\gamma_{12}]\centerdot\chi_k[\overline\zeta_{12}]\centerdot \chi_k[\alpha_{12}]\centerdot\chi_l[\overline\alpha_{01}]
$$
is a simple closed curve and thus, by the Sch\"onflies Theorem, bounds a region $A$ homeomorphic to the unit disk $D^2$.
Similarly, 
$$
\ku K'=\chi_l[\zeta_{10}]\centerdot\chi_l[\gamma_{10}]\centerdot\chi_k[\overline\gamma_{21}]\centerdot\chi_k[\overline\zeta_{21}]\centerdot \chi_k[\alpha_{21}]\centerdot\chi_l[\overline\alpha_{10}]
$$
is a simple closed curve and thus bounds a region $A'$ homeomorphic to the unit disk $D^2$.

Now, as  $\chi_l[\alpha_{01}]\centerdot \chi_k[\overline\alpha_{12}]\subseteq \chi_l[R_{01}]\cup \chi_{k}[R_{12}]$, by condition \ref{g} on $g$, 
$$
g[\chi_l[\alpha_{01}]\centerdot \chi_k[\overline\alpha_{12}]]\subseteq{\rm Int_{M_0}}(\chi_l[R^+_{01}]\cup \chi_{k}[R^+_{12}])
$$ 
and hence intersects $\chi_l[\zeta_{01}]\centerdot \chi_l[\gamma_{01}]\centerdot\chi_k[\overline\gamma_{12}]\centerdot\chi_k[\overline\zeta_{12}]$ only in their common endpoints. Thus, 
$$
\ku L=\chi_l[\zeta_{01}]\centerdot \chi_l[\gamma_{01}]\centerdot\chi_k[\overline\gamma_{12}]\centerdot\chi_k[\overline\zeta_{12}]
\centerdot g[\chi_k[\alpha_{12}]]\centerdot g[\chi_l[\overline\alpha_{01}]]
$$
is a simple closed curve bounding a region $B$ homeomorphic to $D^2$. Similarly, 
$$
\ku L'=\chi_l[\zeta_{10}]\centerdot \chi_l[\gamma_{10}]\centerdot\chi_k[\overline\gamma_{21}]\centerdot\chi_k[\overline\zeta_{21}]
\centerdot g[\chi_k[\alpha_{21}]]\centerdot g[\chi_l[\overline\alpha_{10}]]
$$
bounds a region $B'$ homeomorphic to $D^2$.

We now have two decompositions of $\chi_l[R^+_{01}]\cup \chi_{k}[R^+_{12}]$.
\begin{enumerate}
\item $A\cup \big[\chi_l[R_{01}]\cup \chi_{k}[R_{12}]\big]\cup A'$.
\item $B\cup g[\chi_l[R_{01}]\cup \chi_{k}[R_{12}]]\cup B'$.
\end{enumerate}
Here $A$ and $\chi_l[R_{01}]\cup \chi_{k}[R_{12}]$ overlap along the edge $\chi_l[\alpha_{01}]\centerdot \chi_k[\overline\alpha_{12}]$, $\chi_l[R_{01}]\cup \chi_{k}[R_{12}]$ and $A'$ overlap along $\chi_l[\alpha_{10}]\centerdot \chi_k[\overline\alpha_{21}]$, while $A\cap A'=\tom$. Similarly, $B$ and $g[\chi_l[R_{01}]\cup \chi_{k}[R_{12}]]$ overlap along the edge $g[\chi_l[\alpha_{01}]]\centerdot g[\chi_k[\overline\alpha_{12}]]$, $g[\chi_l[R_{01}]\cup \chi_{k}[R_{12}]]$ and $B'$ overlap along $g[\chi_l[\alpha_{10}]]\centerdot g[\chi_k[\overline\alpha_{21}]]$, while $B\cap B'=\tom$.

We can now define a homeomorphism $\varphi_{lk}: \chi_l[R^+_{01}]\cup \chi_k[R^+_{12}]\til\chi_l[R^+_{01}]\cup \chi_k[R^+_{12}]$, by first setting $\varphi_{lk}=g\inv$ on $g[\chi_l[R_{01}]\cup \chi_k[R_{12}]]$, and then let $\varphi_{lk}$ send $B$ to $A$, while fixing each point of $\chi_l[\zeta_{01}]\centerdot \chi_l[\gamma_{01}]\centerdot\chi_k[\overline\gamma_{12}]\centerdot\chi_k[\overline\zeta_{12}]$ and be $g\inv$ on $g[\chi_l[\alpha_{01}]\centerdot \chi_k[\overline\alpha_{12}]]$. Similarly for $B'$ and $A'$.

This can be done for all pairs of $\chi_l$ and $\chi_k$ with a common edge, and we thus produce homeomorphisms $\varphi_{lk}$ on all of  the regions, similar to $\chi_l[R^+_{01}]\cup \chi_k[R^+_{12}]$, that fix each point of the boundary curve 
$$
\chi_l[\omega_{10}]\centerdot \chi_l[\zeta_{10}]\centerdot \chi_l[\gamma_{01}]\centerdot\chi_k[\overline\gamma_{12}]\centerdot\chi_k[\overline\zeta_{12}]\centerdot \chi_k[\omega_{12}]\centerdot \chi_k[\zeta_{21}]\centerdot\chi_k[\gamma_{21}]\centerdot \chi_l[\overline\gamma_{10}]\centerdot\chi_l[\overline\zeta_{10}].
$$
Pasting all of these $\varphi_{lk}$ together and extending to all of $M$ by setting $\phi={\rm id}$ elsewhere, we obtain  a homeomorphism $\phi\in {\rm Homeo}(M)$ whose support is contained in $\bigcup_{l=1}^m\;\bigcup_{0\leq i<j\leq 2}\chi_l[R^+_{ij}]$, while being the inverse of $g$ on $\bigcup_{l=1}^m\;\bigcup_{0\leq i<j\leq 2}\chi_l[R_{ij}]$. By Lemma \ref{band-support}, $\phi\in W^{20}$, which finishes the proof.
\end{proof}

We are now ready to finish the proof of the Theorem using the preceding sequence of lemmas.
\begin{proof}Let $y_1,\ldots,y_p\in M$ be the vertices of the triangulation and choose for each $i\leq p$ a neighbourhood $U_i$ of $y_i$ homeomorphic to $\R^2$. Find also $0<\eps<\eta_1$ such that $\overline B(y_i,\eps)\subseteq U_i$ for all $i$. By Lemma \ref{patched}, there are $0<\delta_0<\eps$, $z_i\in B(y_i,\eps)$, such that if $g\in{\rm Homeo}(M)$ has support contained in $\bigcup_{i=1}^p\overline B(z_i,\delta_0)$, then $g\in W^{16}$. As $y_i,z_i\in U_i\isom \R^2$, we can, as $W^2$ is dense in $V_{\eta_1}$, find some $h_0\in W^2$ such that $h_0(y_i)\in U'_i\subseteq \overline B(z_i,\delta_0)$, where $U_i'$ is a neighbourhood of $z_i$ homeomorphic to $\R^2$. Therefore, there is some $g_0\in W^{16}$ such that $g_0h_0(y_i)=z_i$. This shows that if $f\in {\rm Homeo}(M)$ has support contained in $U=(g_0h_0)\inv[\bigcup_{i=1}^p]$, then $(g_0h_0)\inv f(g_0h_0)$ has support contained in $\bigcup_{i=1}^p B(z_i,\delta_0)$ and hence belongs to $W^{16}$. So $f$ belongs to $W^{52}$.
We notice that $U$ is an open set containing $y_1,\ldots, y_p$.

Recall now the definition of the paths $\alpha_{ij}, \beta_{ij}$, etc. and also the fact that these paths all depend on the choice of $0<\tau<1$. For a fixed choice of $\tau$, we define the following simple closed curves in $\triangle$
\begin{equation}\begin{split}
\ku F_0^\tau= \beta_{02}\centerdot\theta_{0}\centerdot\overline\beta_{01}\centerdot\overline\alpha_{01}\centerdot\zeta_{01}\centerdot\overline\xi_{0}\centerdot\overline\zeta_{02}\centerdot\alpha_{02},\\
\ku F^\tau_1=\beta_{10}\centerdot\theta_{1}\centerdot\overline\beta_{12}\centerdot\overline\alpha_{12}\centerdot\zeta_{12}\centerdot\overline\xi_{1}\centerdot\overline\zeta_{10}\centerdot\alpha_{10},\\
\ku F_2^\tau=\beta_{21}\centerdot\theta_{2}\centerdot\overline\beta_{20}\centerdot\overline\alpha_{20}\centerdot\zeta_{20}\centerdot\overline\xi_{2}\centerdot\overline\zeta_{21}\centerdot\alpha_{21}.\\
\end{split}\end{equation}
Moreover, we let $F_0^\tau, F_1^\tau, F^\tau_2$ be the closed regions that they enclose. We notice that $F_i^\tau$ converges in the Vietoris topology to $\{v_i\}$ when $\tau\til 0$, and thus for some $\tau>0$, we have for all $i=0,1,2$ and $l=1,\ldots, m$, $\chi_l[F^\tau_i]\subseteq U$. So fix this $\tau$ and denote $F^\tau_i$ by $F_i$. We notice that 
$$
\triangle=(1-\tau)\triangle\cup \bigcup_{0\leq i<j\leq 2}R_{ij} \cup \bigcup_{i=0,1,2}F_i.
$$
By consequence, if $f\in {\rm Homeo}(M)$ is the identity on
$$
\bigcup_{i=1}^m\chi_i[ {(1-\tau)\triangle}]\cup \bigcup_{l=1}^m\;\bigcup_{0\leq i<j\leq 2}\chi_l[R_{ij}],
$$ 
then $f$ has support contained in $\bigcup_{l=1}^m\; \bigcup_{i=0,1,2}\chi_l[F_i]\subseteq U$, and hence $f\in W^{52}$.

Find now a $\nu>0$ as in the statement of Lemma \ref{vertices left}.  Then if $g\in V_\nu$ and $g$ is the identity on $\bigcup_{i=1}^m\chi_i[ {(1-\tau)\triangle}]$, then there is a $\phi\in W^{20}$ such that $\phi\circ g$ is the identity on 
$$
\bigcup_{i=1}^m\chi_i[ {(1-\tau)\triangle}]\cup \bigcup_{l=1}^m\;\bigcup_{0\leq i<j\leq 2}\chi_l[R_{ij}],
$$
and hence belongs to $W^{52}$. But then also $g\in W^{72}$.

Fix $\delta<\frac\nu2$ and find an $\eta>0$ satisfying $\eta<\frac1{1000}$, $\eta<\frac\nu2$, and such that for $i\leq m$ and $x,y\in \triangle$,
$$
d(x,y)<100\eta\til d(\chi_i(x),\chi_i(y))<\delta. 
$$
By Lemma \ref{fixing triangles}, we can find an $0<\alpha<\frac\nu2$ such that for all $h\in V_\alpha$ there is $\psi\in V_\delta\cap W^{20}$ such that for all $i\leq m$,
$$
\psi\circ h|_{\chi_i[(1-2\eta)\triangle]}={\rm id}.
$$
In particular, $\psi\circ h\in V_{\delta}V_{\alpha}\subseteq V_{\delta+\alpha}\subseteq V_\nu$ and is the identity on $\bigcup_{i=1}^m\chi_i[ {(1-\tau)\triangle}]$, whereby $\psi\circ h\in W^{72}$ and $h\in W^{92}$. This shows that $V_\alpha\subseteq W^{92}$ and thus $W^{92}$ contains an open neighbourhood of the identity in ${\rm Homeo}(M)$ and hence we have proved that ${\rm Homeo}(M)$ is Steinhaus, which finishes the proof of the Theorem.
\end{proof}

\begin{flushleft}

Department of mathematics,\\
University of Illinois at Urbana-Champaign,\\
273 Altgeld Hall, MC 382,\\
1409 W. Green Street,\\
Urbana, IL 61801,\\
USA.\\
\texttt{rosendal@math.uiuc.edu}
\end{flushleft}

\end{document}